\documentclass[12pt]{article}
\input epsf
\newcommand{\bC}{{\bf C}}
\newcommand{\bP}{{\bf P}}
\newcommand{\bN}{{\bf N}}
\newcommand{\bR}{{\bf R}}
\newcommand{\bQ}{{\bf Q}}
\newcommand{\cN}{{\cal N}}
\newcommand{\cB}{{\cal B}}

\newcommand{\cQ}{{\cal Q}}
\newcommand{\cH}{{\cal H}}
\def\GCD{\mbox{\rm GCD}}
\def\ord1{\mbox{\rm ord}}
\def\supp{\mbox{\rm supp}\,}
\def\disc{\mbox{\rm disc}}
\def\resultant{\mbox{\rm resultant}}
\def\In{\mbox{\rm in}}            

\newcommand{\Teis}[2]{
   \setlength{\unitlength}{1ex}
   \begin{picture}(2,0)(0,0.4)
      \put(0,1.1){\line(1,0){2}}
      \put(0,0.9){\line(1,0){2}}
      \put(1,1.2){\makebox(0,0)[b]{$\scriptstyle #1$}}
      \put(1,0.8){\makebox(0,0)[t]{$\scriptstyle #2$}}
   \end{picture}}
\newcommand{\Teiss}[2]{
   \setlength{\unitlength}{1ex}
   \begin{picture}(3,2)(0,0.4)
      \put(0,1.1){\line(1,0){3}}
      \put(0,0.9){\line(1,0){3}}
      \put(1.5,1.2){\makebox(0,0)[b]{$\scriptstyle #1$}}
      \put(1.5,0.8){\makebox(0,0)[t]{$\scriptstyle #2$}}
   \end{picture}}
\newcommand{\Teisss}[4]{
   \setlength{\unitlength}{1ex}
   \begin{picture}(#3,3)(0,0.4)
      \put(0,1.15){\line(1,0){#3}}
      \put(0,0.85){\line(1,0){#3}}
      \put(#4,1.3){\makebox(0,0)[b]{$#1$}}
      \put(#4,0.7){\makebox(0,0)[t]{$#2$}}
   \end{picture}}
\newcommand{\Teissss}[2]{
   \setlength{\unitlength}{1ex}
   \begin{picture}(7,3)(0,0.4)
      \put(0,1.15){\line(1,0){7}}
      \put(0,0.85){\line(1,0){7}}
      \put(3.5,1.3){\makebox(0,0)[b]{$#1$}}
      \put(3.5,0.7){\makebox(0,0)[t]{$#2$}}
   \end{picture}}

\newtheorem{thm}{Theorem}[section]
\newtheorem{lm}[thm]{Lemma}
\newtheorem{cor}[thm]{Corollary}
\newtheorem{ex}[thm]{Example}
\newtheorem{prop}[thm]{Proposition}
\newtheorem{rem}[thm]{Remark}
\newtheorem{proper}[thm]{Property}
\title{POLAR INVARIANTS OF PLANE CURVE SINGULARITIES:
INTERSECTION THEORETICAL APPROACH
\footnotetext{
\begin{minipage}[t]{5in}
{\small
2000 {\it Mathematics Subject Classification:\/} Primary 32S55;
Secondary 14H20.
Key words and phrases: plane curve singularity, polar invariant, jacobian
Newton polygon, equisingularity, pencil of plane curve singularities.\\
$^*$The second-named author was partially supported by the KBN grant
No N N201 386634
} 
\end{minipage}
}}

\begin{document}
\author{Janusz Gwo\'{z}dziewicz, Andrzej Lenarcik$^*$
 and Arkadiusz P\l{}oski}

\maketitle

\begin{center}
{\it Dedicated to Professor Pierrette Cassou-Nogu\`{e}s}
\end{center}

\begin{abstract}
\noindent
This article, based on the talk given by one of the authors
at the Pierrettefest in Castro Urdiales in June 2008,
is an overview of a number of recent results
on the polar invariants of plane curve singularities.
\end{abstract}

\section*{Introduction}

The polar invariants (called also polar quotients) of isolated hypersurface
singularities were introduced by B.~Teissier in 1975 to study
equisingularity problems (see \cite{Te1975}, \cite{T3}, \cite{T4}).
They are by definition, the contact orders
between a hypersurface and the branches of its generic polar curve. To every
polar invariant~$q$ of a given isolated hypersurface singularity one
associates in a natural way an integer $m_q>0$ called the multiplicity
of~$q$. Teissier's collection~$\{(q,m_q)\}$ is an analytic invariant
of the singularity. Even more: it is an invariant of the ``c-cos\'ecance''
which is equivalent in the case of plane curve singularities to the
constancy of the local embedded topological type (see \cite{T3}).
The Milnor number, the \L{}ojasiewicz exponent, the ${\cal C}^0$-degree
of sufficiency and other numerical invariants can be computed in terms
of Teissier's collection.

It is well-known (see \cite{T2}, \cite{BK}, \cite{T1}) that the
constancy of the local embedded topological type of plane curves
is equivalent to the usual definitions of equisingularity (see
Preliminaries where the definition of equisingularity in terms of
intersection numbers is given).

M.~Merle \cite{M} computed Teissier's collection for a branch
(irreducible analytic curve) in terms of the semigroup of the
branch. Much earlier a computation of the contacts between an
irreducible curve and the branches of its generic polar curve was
done by Henry J.~S.~Smith \cite{S} but his work fell into oblivion
for long time. R.~Ephraim \cite{Eph} generalized the Smith-Merle
result computing the polar invariants in the case of special
polars and applied his result to the pencil of curves which
appears when studying affine curves with one branch at infinity
(see Sections~4 and~7 of this article).

The case of multi-branched curves turned out much more complicated
and was studied by many authors: Eggers~\cite{E},
Delgado~\cite{D}, Casas-Alvero~\cite{Casas2000}, Garc\'{\i}a
Barroso~\cite{G}, C.~T.~C.~Wall~\cite{Wall1} using algebraic
methods and by L\^e~D.~T., F.~Michel and C.~Weber
in~\cite{LMW1989}, \cite{LMW1991} using topological tools.
L\^e~D.~T. initiated the topological approach to the polar
invariants in~\cite{Le1975}.

C.~T.~C.~Wall gave an account of most results obtained in the
above quoted papers in his book~\cite{Wall2} dealing with
different aspects of the curve singularities.

The goal of this article is to give an overview of a number of recent
results on the polar invariants of plane curve singularities.

In Section~2 we present a refinement of Teissier's invariance theorem
in the case of plane curve singularities. In Section~3 we give an approach
to the polar invariants based on Puiseux series developing the method due to
Kuo and Lu~\cite{KL}.

Section~4 is devoted to the Smith-Merle-Ephraim theorem in the one branch
case and to the irreducibility criterion obtained quite recently by
Garc\'{\i}a Barroso and Gwo\'zdziewicz. (Theorem~4.5 and Corollary~4.6).

In Section~5 we present explicit formulae for the polar invariants in terms
of semigroup of branches and intersection multiplicities due to
Gwo\'zdziewicz and P\l{}oski (Theorem~5.2).
The geometric interpretation of these formulae in terms of the Newton
diagrams associated with many-branched singularity is new.

In Section~6 we recall a result obtained by Lenarcik and
P\l{}oski~(Theorem~6.1) which gives an effective formula for the jacobian
Newton diagram (see Section~2) of a nondegenerate (in the sense of
Kouchnirenko) plane curve singularity. Then, we present in Section~7, some
applications of the polar invariants to pencils of plane curve singularities.

\section{Preliminaries}

In this section we recall some useful notions and results that we need
in this article. The references for this part are \cite{BK}, \cite{Casas2000},
\cite{T1}, \cite{Wall2}.

\subsection{Basic notions}

Let $\bC\{X,Y\}$ be the ring of convergent complex power series in variables $X,Y$.
Let $f\in\bC\{X,Y\}$ be a nonzero power series without constant term.
An {\it analytic curve\/} $f=0$ is defined to be the ideal generated by
$f$ in $\bC\{X,Y\}$. We say that $f=0$ is {\it irreducible (reduced)\/}
if $f\in\bC\{X,Y\}$ is irreducible ($f$ has no multiple factors).
The irreducible curves are also called branches.
If $f=f_1^{m_1}\dots f_r^{m_r}$ with non-associated irreducible
factors $f_i$ then we refer to $f_i=0$ as the branches or components
of $f=0$.

Recall here that for any nonzero power series $f=\sum
c_{\alpha\beta}X^\alpha Y^\beta$ we put
$\ord1\,f=\inf\{\alpha+\beta:\; c_{\alpha\beta}\neq 0\}$ and
$\In\,f=\sum c_{\alpha\beta}X^\alpha Y^\beta$ with summation over
$(\alpha,\beta)$ such that $\alpha+\beta=\ord1\,f$. The initial
form $\In\,f$ of $f$ determines the \textit{tangents to $f=0$}.

For any power series $f,g\in\bC\{X,Y\}$ we define the
{\it intersection number\/} $(f,g)_0$ by putting
$$
  (f,g)_0=\mbox{dim}_{\bC}\bC\{X,Y\}/(f,g)
$$
where $(f,g)$ is the ideal of $\bC\{X,Y\}$ generated by $f$ and $g$.
If $f,g$ are nonzero power series without constant term then
$(f,g)_0<+\infty$ if and only if the curves $f=0$ and $g=0$
have no common branch.

Now suppose that $f=0$ is a branch and consider
$$
  S(f)=\{(f,g)_0:\;g\in\bC\{X,Y\}\mbox{ runs over all series such that
  $f$ does not divide $g$}\}\;.
$$
Clearly $0\in S(f)$ (take $g=1$) and $a,b\in S(f)\Rightarrow a+b\in S(f)$
since the intersection number is additive. We call $S(f)$ the
{\it semigroup of the branch\/} $f=0$. Note that $S(f)=\bN$ if and only if
$\ord1\,f=1$ (we say then that $f=0$ is \textit{regular} or \textit{nonsingular}).

Consider two reduced curves $f=0$ and $g=0$. They are {\it equisingular\/}
if and only if there are factorizations
$f=\prod_{i=1}^r f_i$ and $g=\prod_{i=1}^r g_i$ with the same number $r>0$
of irreducible factors $f_i$ and $g_i$ such that
\begin{itemize}
\item  $S(f_i)=S(g_i)$ for all $i=1,\dots,r$,
\item  $(f_i,f_j)_0=(g_i,g_j)_0$ for $i,j=1,\dots,r$.
\end{itemize}
The bijection $f_i\mapsto g_i$ will be called
{\it equisingularity bijection\/}.
In particular two branches are equisingular if and only if they have
the same semigroup. A function defined on the set of reduced curves is an
{\it invariant\/} if it is constant on equisingular curves.
The multiplicity $\ord1\,f$, the number of branches $r(f)$ and the
number of tangents $t(f)$ of $f=0$ are invariants.

For any analytic curve $f=0$ we consider the {\it Milnor number\/}
$\mu_0(f)=(\partial f/\partial X,\partial f /\partial Y)_0$. One has $\mu_0(f)<+\infty$ if and
only if the curve $f=0$ is reduced. Let us recall the following two
properties:
\begin{itemize}
\item  if $f=0$ is a branch then $\mu_0(f)$ is the smallest integer
$c\geq 0$ such that all integers greater than or equal to $c$
belong to $S(f)$,
\item if $f=f_1\dots f_r$ with pairwise different
irreducible $f_i$ then
$$
  \mu_0(f)+r-1=\sum_{i=1}^r\mu_0(f_i)+2\sum_{1\leq i<j\leq r}(f_i,f_j)_0\;.
$$
\end{itemize}
Thus the Milnor number is an invariant. A simple proof of the above
properties is given in~\cite{P1}.

\subsection{Newton diagrams after~\cite{T2}}

Let $\bR_+=\{x\in\bR:\;x\geq 0\}$. The Newton diagrams are some
convex subsets of $\bR_+^2$. Let $E\subset\bN^2$ and let us denote
by $\Delta(E)$ the convex hull of the set $E+\bR_+^2$. The subset
$\Delta$ of $\bR_+^2$ is a {\it Newton diagram\/} (or
\textit{polygon}) if there is a set $E\subset\bN^2$ such that
$\Delta=\Delta(E)$. The smallest set $E_0\subset\bN^2$ such that
$\Delta=\Delta(E_0)$ is called the \textit{set of vertices} of the
Newton diagram $\Delta$. It is always finite and we can write
$E_0=\{v_0,v_1,\dots,v_m\}$ where $v_i=(\alpha_i,\beta_i)$ and
$\alpha_{i-1}<\alpha_i$, $\beta_{i-1}>\beta_i$ for all
$i=1,\dots,m$. In particular the Newton diagram $\Delta$ with one
vertex $v=(\alpha,\beta)$ is the quadrant
$(\alpha,\beta)+\bR_+^2$.

According to Teissier for $k,l>0$ we denote $\{\Teis{k}{l}\}$
the Newton diagram with vertices $(0,l)$ and $(k,0)$. We put also
$\{\Teis{k}{\infty}\}=(k,0)+\bR_+^2$ and
$\{\Teis{\infty}{l}\}=(0,l)+\bR_+^2$ and call any diagram of the form
$\{\Teis{k}{l}\}$ an {\it elementary Newton diagram\/}. For any subsets
$\Delta,\Delta'\subset\bR_+^2$ we consider the {\it Minkowski sum\/}
$\Delta+\Delta'=\{u+v:\;u\in\Delta\mbox{ and }v\in\Delta'\}$.
One checks the following
\begin{proper}
The Newton diagrams form the semigroup
with respect to the Minkowski sum. The elementary Newton diagrams generate
the semigroup of the Newton diagrams.
\end{proper}
For any Newton diagram $\Delta$ we consider the set $\cN(\Delta)$
of the compact faces of the boundary of $\Delta$. If $\Delta$ has
vertices $v_0,v_1,\dots,v_m$ then
$\cN(\Delta)=\{[v_{i-1},v_i]:\;i=1,\dots,m\}$. For any segment
$S\in\cN(\Delta)$ we denote by $|S|_1$ and $|S|_2$ the lengths of
the projections of $S$ on the horizontal and vertical axes. We
call $|S|_1/|S|_2$ the \textit{inclination of $S$}. If $\Delta$
intersects both axes then
$\Delta=\sum_S\left\{\Teiss{|S|_1}{|S|_2}\right\}$ (summation over
all $S\in\cN(\Delta)$) and this representation is unique.

Now, let $f=\sum c_{\alpha\beta}X^\alpha Y^\beta$ be a power series.
We put $\supp\,f=\{(\alpha,\beta)\in\bN^2:\,c_{\alpha\beta}\neq 0\}$,
$\Delta_{X,Y}(f)=\Delta(\supp\,f)$ and $\cN_f=\cN(\Delta(f))$.
We call $\Delta_{X,Y}(f)$ the \textit{Newton diagram (or polygon) of the power
series $f$}. Let $n>0$ be an integer. Let $f=f(X,Y)$ be a power series
$Y$-regular of order $n$, i.e. such that $\ord1\,f(0,Y)=n$.
Let $\bC\{X\}^*=\bigcup_{p\geq 1}\bC\{X^{1/p}\}$ be the ring of Puiseux
series. We have the Newton-Puiseux factorization
$$
  f(X,Y)=U(X,Y)\prod_{i=1}^n(Y-\alpha_i(X)),\quad U(X,Y)\mbox{ is a unit in }\bC\{X,Y\}
$$
where $\alpha_i(X)\in\bC\{X\}^*$ for $i=1,\dots,n$.
\begin{thm} {\rm(Newton-Puiseux Theorem)}\\
For every $q\in\bQ\cup \{\infty\}$ let $m_q$ be the number of
roots $\alpha_i(X)$ such that $\ord1\,\alpha_i(X)=q$. Then $m_qq$
(by convention $0\cdot\infty=0$) is an integer or $\infty$ and
$$
\Delta_{X,Y}(f)=\sum_{q}\left\{
         \Teisss{m_q q}{m_q}{5}{2.5}
         \right\}\;.
$$
\end{thm}
\subsection{Nondegeneracy}

Now, let $f=\sum c_{\alpha\beta}X^\alpha Y^\beta$ be a power series.
For any segment $S\in\cN(f)$ we put
$\In(f,S)=\sum c_{\alpha\beta}X^\alpha Y^\beta$ where $(\alpha,\beta)\in S$.

According to~\cite{K} , the series $f$ is \textit{nondegenerate} if for every
$S\in\cN(f)$ the polynomial $\In(f,S)$ has no critical points in the set
$\bC^*\times\bC^*$, where $\bC^*=\bC\setminus\{0\}$.
A lot of applications of the Newton diagrams are based on the following
\begin{thm}\label{thm03} {\rm (\cite{GLP0}, \cite{L3}).}
Suppose that $f,g\in\bC\{X,Y\}$
are reduced power series such that $\Delta(f)=\Delta(g)$. Then
\begin{itemize}
\item[\rm(i)] if $f$ and $g$ are nondegenerate then the curves $f=0$ and
$g=0$ are equisingular,
\item[\rm(ii)] if $f$ is nondegenerate but $g$ is degenerate then $f=0$
and $g=0$ are not equisingular.
\end{itemize}
\end{thm}
Let $\Delta\subset\bR_+^2$ be a Newton diagram. It is easy to
check that $\Delta=\Delta(f)$ for a reduced nondegenerate power
series $f$ if and only if the distances from $\Delta$ to the axes
are~${}\leq1$. We call such diagrams \textit{nearly convenient}.
Every Newton diagram which intersects both axes (convenient in the
sense of Kouchnirenko) is nearly convenient. If $\Delta$ is nearly
convenient then the reduced nondegenerate power series~$f$ such
that~$\Delta=\Delta(f)$ form an open dense subset in the space of
coefficients.

Let us consider an invariant $I$ of equisingularity.  For every
nearly convenient  Newton diagram $\Delta$ we put
$I(\Delta)=I(\Delta(f))$ where $f$ is a nondegenerate reduced
power series. According to the theorem quoted above $I(\Delta)$ is
defined correctly (does not depend on $f$). There is a natural
problem: calculate $I(\Delta)$ effectively in terms of $\Delta$.
The most known result of this kind is due to
Kouchnirenko~\cite{K}.

To formulate it let us consider for every nearly convenient Newton
diagram~$\Delta$ a convex subset $\tilde\Delta$ of~$\bR_+^2$
defined to be the intersection of all half-planes
containing~$\Delta$ whose boundary is the line extending a
face~$S\in{\cN}(\Delta)$ with $\bR_+^2$. If
${\cN}(\Delta)=\emptyset$ then by convention
$\tilde\Delta=\bR_+^2$. Let~$(a,0)$ (resp.~$(0,b)$) be the point
of $\tilde\Delta\cap\{\beta=0\}$
(resp.~$\tilde\Delta\cap\{\alpha=0\}$) closest to the origin. Let
us put $\mu(\Delta)= 2\cdot\mbox{area of
}(\bR_+^2\setminus\tilde\Delta)-a-b+1$. Then we have

\begin{thm}\label{thK}{\rm (see \cite{K}, \cite{GLP0}).}
For any power series $f$: $\mu_0(f)\geq\mu(\Delta(f))$.
The equality holds if and only if $f$ is nondegenerate.
\end{thm}

\vspace{1ex}\noindent Note that Kouchnirenko proved a much more general result
concerning isolated singularities in $n$ dimensions. In the case $n=2$
the result is more precise: the equality $\mu_0(f)=\mu(\Delta(f))$ holds
if and only if $f$ is nondegenerate and we do not need the assumption
``$f$ is convenient''.

Theorem 1.3 can be easily deduced from the famous
$\mu$-constant theorem~\cite{LeR} and Kouchnirenko's result.
One can give also a direct, elementary proof~\cite{L3}.
Let us end this section with
\begin{ex}
{\rm Let $f(X,Y)=\sum c_{\alpha\beta}X^\alpha Y^\beta$
($\frac{\alpha}{w_1}+\frac{\beta}{w_2}=1$ where $w_1,w_2\geq 2$
are rational numbers) be a weighted homogeneous polynomial of
order$\mbox{}>1$. Then $\bR_+^2\setminus\tilde\Delta(f)$ is the
triangle with sides $\alpha=0$, $\beta=0$ and
$\alpha/w_1+\beta/w_2=1$.  If $f$ is nondegenerate then by
Theorem~\ref{thK} $\mu_0(f)=\mu(\Delta(f))=(w_1-1)(w_2-1)$ (the
Milnor-Orlik formula).
} 
\end{ex}

\section{The jacobian Newton polygon}

The following lemma is well-known (see, for example~\cite{D1991}
or~\cite{P3}).
\begin{lm}\label{lm11}
Let $f,g\in\bC\{X,Y\}$ be two power series without constant term.
Let $J(f,g)=({\partial f}/{\partial X})({\partial g}/{\partial Y})-
({\partial f}/{\partial Y})({\partial g}/{\partial X})$ be the Jacobian of the pair $(f,g)$.
Then
$$
  (f,J(f,g))_0=\mu_0(f)+(f,g)_0-1\;.
$$
The right side of the above equality is finite if and only if the left
is too.
\end{lm}
Assume that $l=0$ is a regular curve.
Let $f=0$ be a reduced curve such that $J(f,l)(0,0)=0$.
If $l=0$ is not a branch of $f=0$
then we call $J(f,l)=0$ the \textit{polar curve of $f=0$ relative to
$l=0$}. It depends on the power series $f$ and $l$.

If $l=bX-aY$ is a nonzero linear form then
$$
J(f,l)=a({\partial f}/{\partial X})+b({\partial f}/{\partial Y})
$$
and we speak about the \textit{polar curve relative to the
direction $(a:b)\in\bP^1(\bC)$}. Using Lemma~\ref{lm11} it is easy
to check the following two properties. We assume $J(f,l)(0,0)=0$.
\begin{proper}\label{proper12}
The regular curve $l=0$ is not a branch of the curve $f=0$ if and only if
$l=0$ is not a branch of the polar curve $J(f,l)=0$.
\end{proper}
Recall that two curves are \textit{transverse} if they have no
common tangent.
\begin{proper}\label{proper13}
If the curves $l=0$ and $f=0$ are transverse then the curves
$l=0$ and $J(f,l)=0$ are transverse, too.
\end{proper}
In the sequel we assume that $f=0$ is a reduced curve and that the
regular curve $l=0$ is not a branch of $f=0$.

Recall that $J(f,l)(0,0)=0$ and let $J(f,l)=h_1\cdots h_s$ be the
decomposition of $J(f,l)$ into irreducible factors. Then the
rational numbers
$$
  \frac{(f,h_j)_0}{(l,h_j)_0},\;j=1,\dots,s
$$
are called the {\it polar invariants\/} of $f=0$ relative to $l=0$.
Let $Q(f,l)$ be the set of polar invariants. If $J(f,l)(0)\neq 0$
then we put $Q(f,l)=\emptyset$. For every polar invariant
$q\in Q(f,l)$ we put
$$
  A_q=\{j\in[1,s]:\;\frac{(f,h_j)_0}{(l,h_j)_0}=q\}
$$
and
$$
  J_q=\prod_{j\in A_q}h_j\;.
$$
Thus
$$
J(f,l)=\prod_q J_q\mbox{\quad and\quad }
\frac{(f,J_q)_0}{(l,J_q)_0}=q \quad \mbox{for } q\in Q(f,l)\;.
$$
We call $m_q=(l,J_q)_0$ the {\it multiplicity\/} of the polar invariant $q$.
Using Lemma~\ref{lm11} we check
\begin{proper}\label{proper14}
$\displaystyle
\quad\sum_q m_q=(f,l)_0-1,\quad
\sum_q m_q q=\mu_0(f)+(f,l)_0-1$ where the summation is over
all $q\in Q(f,l)$.
\end{proper}
Let $\eta_0(f,l)=\sup\,Q(f,l)$ be the maximal polar invariant
($\eta_0(f,l)=-\infty$ if $J(f,l)(0,0)\neq 0$).
Property~\ref{proper14} implies
\begin{proper}\label{proper15}
Suppose that $(f,l)_0>1$. Then
$$
\frac{\mu_0(f)}{(f,l)_0-1}+1\leq\eta_0(f,l)\leq\mu_0(f)+1\;.
$$
Moreover
$$
\eta_0(f,l)=\frac{\mu_0(f)}{(f,l)_0-1}+1
$$
if and only if there exists exactly one polar invariant of $f=0$
relative to $l=0$.
\end{proper}
From the above property it follows that a regular plane curve $f=0$
has exactly one polar invariant, equal to 1 relative to any
nontransverse regular curve $l=0$.
In the sequel we assume that $f=0$ is a singular reduced curve.

Following Teissier~\cite{T4} we define the \textit{jacobian Newton polygon} by putting
$$
\cQ(f,l)=\sum_{j=1}^s\left\{
         \Teissss{(f,h_j)_0}{(l,h_j)_0}
         \right\}\;.
$$
It is easy to see that
$$
\cQ(f,l)=\sum_{q}\left\{
         \Teisss{m_q q}{m_q}{5}{2.5}
         \right\}\;.
$$
\begin{proper}\label{proper16}
The jacobian Newton polygon intersects
the axes at points $(0,(f,l)_0-1)$ and $(\mu_0(f)+(f,l)_0-1,0)$.
All faces of $\cQ(f,l)$ have inclinations
strictly greater than $1$.
\end{proper}
The above property follows from Property~\ref{proper14} and from the
following formula
$$
  (f,h_j)_0=\inf\left\{
  \left(\frac{\partial f}{\partial X},h_j\right)_0,
  \left(\frac{\partial f}{\partial Y},h_j\right)_0
  \right\}+(l,h_j)_0\mbox{ for }j=1,\dots,s\;.
$$
$$
\epsfbox{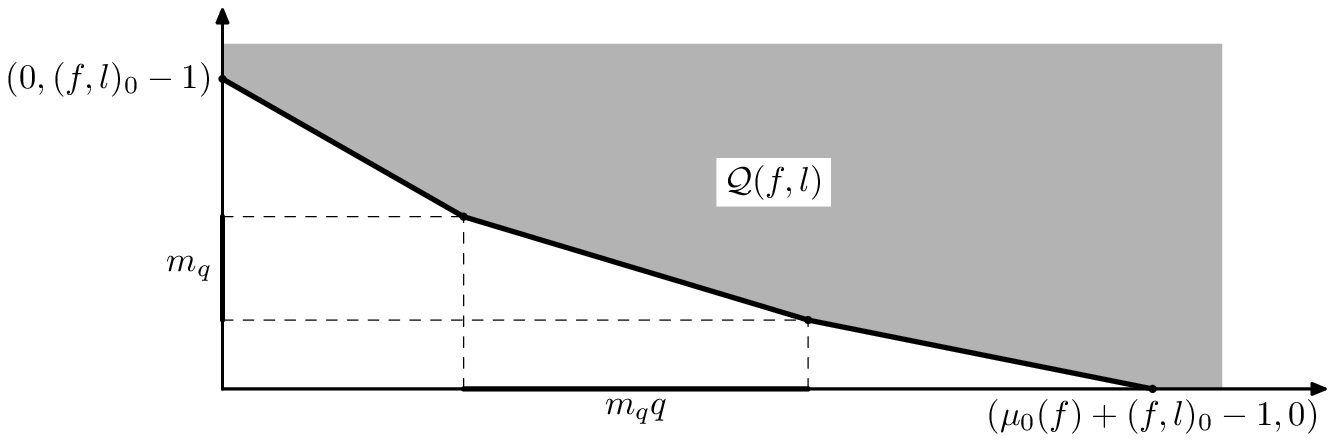}
$$
\begin{rem}\label{rem17}
{\rm If $f=0$ and $l=0$ are transverse then the polar invariants are of the
form $(f,h_j)_0/\ord1\,h_j$, $j=1,\dots,s$ (see Property~\ref{proper13}).
The jacobian Newton polygon joins the points $(0,\ord1\,f-1)$ and
$(\mu_0(f)+\ord1\,f-1,0)$.
Its faces have inclinations greater than or equal to $\ord1\,f$.
One checks that $\ord1\,f$ is the polar invariant if and only if
the number of tangents $t(f)$ is strictly greater than 1.
Then $t(f)-1$ is the multiplicity of $\ord1\,f$ (see~\cite{LMP}).
} 
\end{rem}

A local isomorphism $\Phi$ is a pair of power series without
constant term such that $\mbox{Jac}\,\Phi(0,0)\neq 0$. The
jacobian Newton polygon $\cQ(f,l)$ is an analytic invariant of the
pair $(f,l)$:
\begin{proper}\label{proper18}
Let $\Phi$ be a local isomorphism. Then
$\cQ(f\circ\Phi,l\circ\Phi)=\cQ(f,l)$.
\end{proper}
Let $f=0$ and $f'=0$ be reduced singular curves and let $l=0$ and $l'=0$
be regular branches such that $l=0$ (resp. $l'=0$) is not a component
of $f=0$ (resp. $f'=0$). We will say that the {\it pairs $f=0,l=0$ and
$f'=0,l'=0$ are equisingular\/} if there is an equisingularity bijection
of the set of branches $f_i=0$ of $f=0$ and $f'_i=0$ of $f'=0$
such that $(f_i,l)_0=(f'_i,l')_0$ for all $i=1,\dots,r$.
The following result is a refinement of Teissier's theorem on invariance
of the jacobian Newton polygon~\cite{T3} in the case of plane curve singularities.
\begin{thm}
Suppose that the pairs $f=0,l=0$ and $f'=0,l'=0$ are equisingular. Then
$$
  \cQ(f,l)=\cQ(f',l')\;.
$$
\end{thm}
The proof of the above theorem may be given by purely intersection
theoretical methods (see~\cite{GLP}) based on the Kuo and Lu
approach (\cite{KL} and Section~3 of this paper).

Now let us note
\begin{cor}\label{cor110}
If $f=0$ is a reduced singular curve and $l=0,l'=0$ are nonsingular branches
transverse to $f=0$ then $\cQ(f,l)=\cQ(f,l')$.
\end{cor}
We write $\cQ(f)=\cQ(f,l)$ provided that $f=0$ and $l=0$ are transverse
and call $\cQ(f)$ the \textit{jacobian Newton polygon of $f=0$}.
\begin{cor}\label{cor111}
Suppose that the reduced singular curves $f=0$ and $f'=0$ are equisingular.
Then $\cQ(f)=\cQ(f')$.
\end{cor}
From the last corollary it follows that the Milnor number $\mu_0(f)$
and the maximal polar invariant $\eta_0(f)=\max\,Q(f)$ are invariants.
\begin{ex}
{\rm Let $f=(Y^3-X^5)^2-9X^{11}$ and $l=X$. Then $(f,l)_0=\ord1\,f=6$
i.e. $f=0$ and $l=0$ are transverse. We get
$J(f,l)=(\partial f/\partial Y)=6(Y^3-X^5)Y^2$ and
$$
  \cQ(f)=\cQ(f,l)=
\left\{\Teisss{(f,Y)_0}{1}{7}{3.5}\right\}+
\left\{\Teisss{(f,Y)_0}{1}{7}{3.5}\right\}+
\left\{\Teisss{(f,Y^3-X^5)_0}{3}{13}{6.5}\right\}=
\left\{\Teisss{20}{2}{3}{1.5}\right\}+
\left\{\Teisss{33}{3}{3}{1.5}\right\}\;.
$$
} 
\end{ex}

The computations of the jacobian Newton polygons in the next two examples
were done using Theorem~6.1.
\begin{ex}
{\rm ([Len2008]) Let $f=Y^9+X^2Y^2+X^9$ and $g=Y^5+XY^4+X^9$. Then
$$
  \cQ(f)=\cQ(g)=
\left\{\Teisss{5}{1}{2}{1}\right\}+
\left\{\Teisss{27}{3}{3}{1.5}\right\}
$$
but the curves $f=0$ and $g=0$ are not equisingular.
The curve $f=0$ has 3 branches while $g=0$ has 5.
} 
\end{ex}
\begin{ex}\label{ex114}
{\rm Let $f=X^3Y^3+X^2Y^4+X^8+Y^7$ and $g=X^4Y^2+X^8+Y^7$. Then
$$
  \cQ(f)=
\left\{\Teisss{6\cdot 2}{2}{4}{2}\right\}+
\left\{\Teisss{7\cdot 1}{1}{4}{2}\right\}+
\left\{\Teisss{8\cdot 2}{2}{4}{2}\right\}
$$
and
$$
  \cQ(g)=
\left\{\Teisss{6\cdot 1}{1}{4}{2}\right\}+
\left\{\Teisss{7\cdot 3}{3}{4}{2}\right\}+
\left\{\Teisss{8\cdot 1}{1}{4}{2}\right\}\;.
$$
We get $\ord1\,f=\ord1\,g=6$ and $\mu_0(f)=\mu_0(g)=30$. The Newton
polygons $\cQ(f)$ and $\cQ(g)$ have the same inclinations $6$, $7$, $8$ and
join the same points $(0,5)$ and $(35,0)$ but $\cQ(f)\neq\cQ(g)$.
} 
\end{ex}
The following simple proposition gives an effective way of computing
the jacobian Newton polygon of the pair
$f(X,Y)=Y^n+a_1(X)Y^{n-1}+\dots+a_n(X)$ (a distinguished polynomial
of degree $n>1$) and $l(X,Y)=X$ by performing the rational operations
on the coefficients $a_1(X),\dots,a_n(X)$.
It illustrates the leading principle of Teissier's lectures~\cite{T2}.

\begin{prop}\label{prop115}
Suppose that $f(X,Y)$ is an $Y$-distinguished polynomial of degree $n>1$
without multiple factors. Let $T$ be a new variable and consider the
discriminant $D(X,T)=\disc_Y(f(X,Y)-T)$. Then
$\cQ(f,X)=\Delta_{X,T}(D)$ {\rm(}the Newton polygon of the discriminant
$D(X,T)$ in coordinates $X,T${\rm)}.
\end{prop}
Proof. Let $\beta_1(X),\dots,\beta_{n-1}(X)\in\bC\{X\}^*$ be the Puiseux roots
of equation $(\partial f/\partial Y)(X,Y)=0$. It is easy
to see that $\ord1\,f(X,\beta_1(X))$,\dots, $\ord1\,f(X,\beta_{n-1}(X))$ is the
sequence of polar invariants of $f=0$ relative to $X$ appearing with their
multiplicities (if $h(X,Y)=0$ is the minimal analytic equation of the
series $\beta(X)\in\bC\{X\}^*$ then $\ord1\,f(X,\beta(X))=(f,h)_0/(X,h)_0$).
On the other hand
\begin{eqnarray*}
D(X,T) & = & \disc_Y(f(X,Y)-T)=\resultant_Y(f(X,Y)-T,\frac{\partial f}{\partial Y}(X,Y))\\
& = & \pm\prod_{j=1}^{n-1}(T-f(X,\beta_j(X)))\;.
\end{eqnarray*}
We apply the Newton-Puiseux Theorem (see Preliminaries) to
$D(X,T)\in\bC\{X,T\}$.
\begin{ex}\label{e116}
{ \rm Let $f(X,Y)=(Y^2-X^3)^2-X^5Y$. Then $f=0$ and $X=0$ are transverse.
We have $D(X,T)=-256T^3+256X^6T^2+288X^{13}T-27X^{20}-256X^{19}$ and
} 
$$
\cQ(f)=\cQ(f,X)=\Delta_{X,T}(D)=
\left\{\Teisss{6}{1}{2}{1}\right\} +
\left\{\Teisss{13}{2}{3}{1.5}\right\}\; .
$$
\end{ex}

\section{Polar invariants and Puiseux series}

The following lemma due to Kuo and Lu (see~\cite{KL}, Lemma~3.3)
is crucial for the approach to the polar invariants based on
Puiseux series (see~\cite{E}, \cite{GP3}, \cite{Wall1}).
\begin{lm}\label{lm21}
{\rm(the Kuo and Lu lemma)}\\
Let $f=f(X,Y)\in\bC\{X,Y\}$ be a $Y$-regular power series of order $n>1$
and let $\alpha_1=\alpha_1(X),\dots,\alpha_n=\alpha_n(X)$ be the Puiseux
roots of the equation $f(X,Y)=0$. If
$\beta_1=\beta_1(X),\dots,\beta_{n-1}=\beta_{n-1}(X)$ are the Puiseux roots
of the equation $(\partial f/\partial Y)(X,Y)=0$ then for each $k\in\{1,\dots,n\}$
and for each $r>0$
$$
  \#\{i:\,\ord1(\alpha_i-\alpha_k)=r\}=
  \#\{i:\,\ord1(\beta_i-\alpha_k)=r\}
$$
\end{lm}
A short proof of the above lemma is given in~\cite{GP1} (see also~\cite{GP3}).
\begin{rem}\label{rem22}
{\rm In~\cite{KL} the following property is stated:
\begin{itemize}
\item[(*)] {\it for given $\alpha_i,\beta_k$ there exists an
$\alpha_j$ such that
$\ord1(\beta_k-\alpha_i)=\ord1(\beta_k-\alpha_j)
=\ord1(\alpha_i-\alpha_j)$.\/}
\end{itemize}
To show that (*) does not hold take $f(X,Y)=Y(Y-X)(Y-X^2)$.
Then $\alpha_1=0$, $\alpha_2=X$, $\alpha_3=X^2$
and $\beta_1=\frac 23X+\dots$, $\beta_2=\frac 12X^2+\dots$.
For $\alpha_2,\beta_2$ does not exist $\alpha_j$ such that
$\ord1(\beta_2-\alpha_2)=\ord1(\beta_2-\alpha_j)
=\ord1(\alpha_2-\alpha_j)$.

Note also that property (*) does not hold under the assumption added in~\cite{G}
that $f(X,0)f(0,Y)\neq 0$.
To get an example it suffices to replace the series
$f(X,Y)$ considered above by the series $f(X,Y-X)$.
} 
\end{rem}

The set of all Puiseux series~$\bC\{X\}^*$ is an ultrametric space
with the order of contact
$O(\varphi,\psi)=\ord1(\varphi(X)-\psi(X))$. That is for any
$\varphi,\psi,\chi\in\bC\{X\}^*$:
\begin{itemize}
\item[(i)] $O(\varphi,\psi)=+\infty$ if and only if $\varphi=\psi$,
\item[(ii)] $O(\varphi,\psi)=O(\psi,\varphi)$,
\item[(iii)] $O(\varphi,\psi)\geq\inf\{O(\varphi,\chi),O(\psi,\chi)\}$.
\end{itemize}
Let $Z\subset\bC\{X\}^*$ be a nonempty finite subset of $\bC\{X\}^*$.
A \textit{ball} in $Z$ is a subset $B\subset Z$ for which there are
$\varphi,\psi\in Z$ such that $\alpha\in B$ if and only if
$O(\alpha,\varphi)\geq O(\varphi,\psi)$. We will write
$B=B(\varphi,O(\varphi,\psi))$. For each ball $B$ in $Z$ we define
the diameter $h(B)=\inf\{O(\alpha,\beta):\,\alpha,\beta\in B\}$.
Note that if $B=B(\varphi,O(\varphi,\psi))$ then
$h(B)=O(\varphi,\psi)$. Let $\cB(Z)$ be the set of balls in $Z$.
The ordered set $(\cB(Z),\leq)$ where $B\leq B'$ if and only if
$B\supset B'$ will be called the \textit{tree over~$Z$}. If $B\leq B'$ with
$B\neq B'$ and there is no other ball between $B$ and $B'$ then we
call $B'$ a {\it successor\/} of $B$. If $h(B)<+\infty$ i.e. if $B$ does
not reduce to a one-point set then $B$ has a finite number $t(B)$ of
successors. One has $t(B)\geq 2$.

Let $f=f(X,Y)\in\bC\{X,Y\}$ be a $Y$-regular power series of order
$n=\ord1\,f(0,Y)\geq 1$. Assume that $f$ has no multiple factors
and let
$$
Z_f=\{\alpha=\alpha(X)\in\bC\{X\}^*:\,\ord1\,\alpha(X)>0
\mbox{ and }f(X,\alpha(X))=0\}\;.
$$
Thus $\# Z_f=n$. The tree over $Z_f$ will be denoted $T(f)$ and called
the \textit{Kuo-Lu tree model of $f$} (see~\cite{KL} where the balls are called bars
and $h(B)$ is called height of $B$).
\begin{ex}\label{ex23}{\rm
(see~\cite{IKK}) Let $f(X,Y)=(Y-X^2)(Y^2-X^3)(Y^2-X^5)$. Here
$\alpha_1=X^2$, $\alpha_2=X^{3/2}$, $\alpha_3=-X^{3/2}$,
$\alpha_4=X^{5/2}$, $\alpha_5=-X^{5/2}$ are the roots of $f(X,Y)=0$.
Thus $Z_f=\{\alpha_1,\dots,\alpha_5\}$ and
$O(Z_f\times Z_f)=\{3/2,2,5/2,+\infty\}$. It is easy to check that
$T(f)=\{B_0,B_1,B_2,\{\alpha_1\},\dots,\{\alpha_5\}\}$ where
$B_0=Z_f$, $B_1=\{\alpha_1,\alpha_4,\alpha_5\}$,
$B_2=\{\alpha_4,\alpha_5\}$. The successors of $B_0$ are
$B_1$, $\{\alpha_2\}$, $\{\alpha_3\}$, the successors of $B_1$ are
$\{\alpha_1\}$ and $B_2$ and the successors of $B_2$ are
$\{\alpha_4\}$ and $\{\alpha_5\}$. Thus we have $t(B_0)=3$,
$t(B_1)=2$, $t(B_2)=2$.
We can represent the tree $T(f)$ in the following figure
$$
\epsfbox{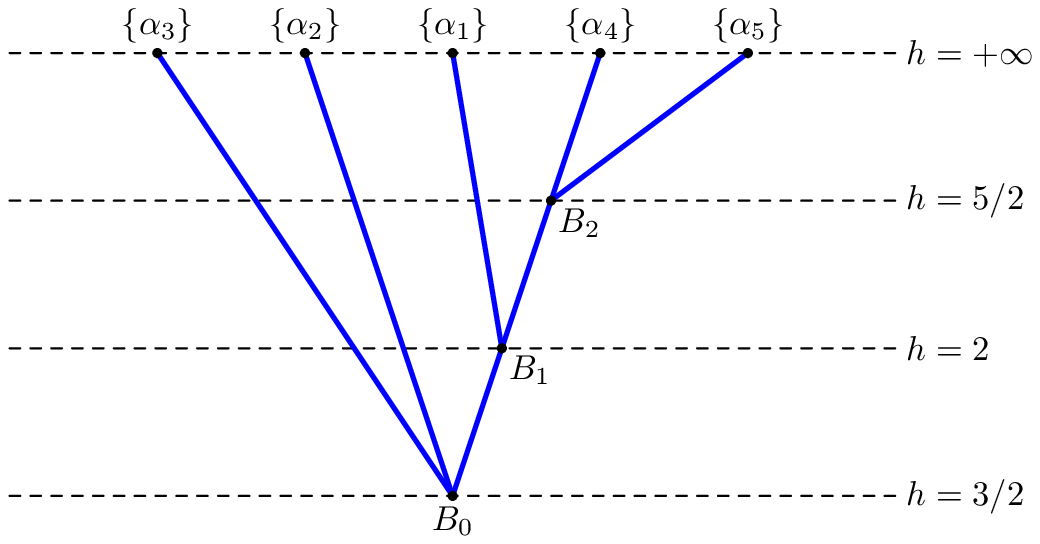}
$$
The balls are represented by points situated on different
levels corresponding to the heights $h\in O(Z_f\times Z_f)$.
We join every ball by continuous lines with its successors.
} 
\end{ex}
For each $\alpha\in Z_f$ and for each ball $B\in T(f)$ we put
$O(\alpha,B)=\sup\{O(\alpha,\varphi):\,\varphi\in B\}$.
Let $T(f)'=\{B\in T(f):\,h(B)<+\infty\}$ and put
$$
  q(B)=\sum_{\alpha\in Z_f}\inf\{O(\alpha,B),\,h(B)\}.
$$
Note that $O(\alpha,B)=O(\alpha,\varphi)$ for any $\varphi\in B$
provided that $O(\alpha,B)<h(B)$.
\begin{thm}
Let $f=f(X,Y)\in\bC\{X,Y\}$, $n=\ord1\,f(0,Y)>1$ be a power series without
multiple factors. Then
\begin{itemize}
\item[\rm(i)] $q\in Q(f,X)$ if and only if $q=q(B)$ for a ball $B\in T(f)'$,
\item[\rm(ii)] $m_q=\sum_B(t(B)-1)$ where summation is over all $B\in T(f)'$
such that $q=q(B)$.
\end{itemize}
\end{thm}
The above quoted theorem is implicit in~\cite{KL}. Part (i) was proved
in~\cite{GP3}. A short proof of (i) and (ii) is given in~\cite{GG}.
\begin{ex}{\rm
Let us calculate $\cQ(f,X)$ for $f=(Y-X^2)(Y^2-X^3)(Y^2-X^5)$.
Using the notation from Example~\ref{ex23} we get
$q(B_0)=(\# Z_f)h(B_0)=5\cdot(3/2)=15/2$,
$q(B_1)=O(\alpha_2,B_1)+O(\alpha_3,B_1)+(\#B_1)h(B_1)=(3/2)+(3/2)+3\cdot 2=9$,
$q(B_2)=O(\alpha_1,B_2)+O(\alpha_2,B_2)+O(\alpha_3,B_2)+(\#B_2)h(B_2)=
2+(3/2)+(3/2)+2\cdot(5/2)=10$.
Consequently, we get
\begin{eqnarray*}
  \cQ(f) = \cQ(f,X) & = &
\left\{\Teisss{(15/2)(3-1)}{3-1}{13}{6.5}\right\}+
\left\{\Teisss{9(2-1)}{2-1}{8}{4}\right\}+
\left\{\Teisss{10(2-1)}{2-1}{9}{4.5}\right\}\\
& = & \left\{\Teisss{15}{2}{4}{2}\right\}+
\left\{\Teisss{9}{1}{3}{1.5}\right\}+
\left\{\Teisss{10}{1}{4}{2}\right\}\;.
\end{eqnarray*}
} 
\end{ex}

\begin{rem}\label{rem26}
{ \rm In~\cite{L2} the polar invariants and their multiplicities
are computed by using the Newton algorithm. }
\end{rem}

\section{The case of one branch}

Let $f=0$ be a singular branch. For any regular curve $l=0$ the
semigroup $S(f)$ has the $(f,l)_0$-minimal system of generators
$\bar{b}_0,\bar{b}_1,\dots,\bar{b}_h$ defined by conditions
\begin{itemize}
\item[\rm(i)] $\bar{b}_0=(f,l)_0$,
\item[\rm(ii)] $\bar{b}_k=\min(S(f)\setminus
(\bN\,\bar{b}_0+\dots+\bN\,\bar{b}_{k-1}))$,
\item[\rm(iii)] $S(f)=\bN\,\bar{b}_0+\dots+\bN\,\bar{b}_{h}$.
\end{itemize}
We will write $\langle \bar{b}_0,\dots,\bar{b}_h\rangle$ instead of
$\bN\,\bar{b}_0+\dots+\bN\,\bar{b}_{h}$.
If $f=0$ and $l=0$ are transverse then $(f,l)_0=\ord1\,f$ and the
corresponding system of $(f,l)_0$-minimal generators will be denoted
$\bar\beta_0,\bar\beta_1,\dots,\bar\beta_g$.
Here $\bar\beta_0=\min(S(f)\setminus\{0\})$.
Let $n_1,\dots,n_h$ be the integers defined to be
$$
  n_k=\frac{\GCD(\bar{b}_0,\dots,\bar{b}_{k-1})}
           {\GCD(\bar{b}_0,\dots,\bar{b}_k)}\quad
           \mbox{ for }k=1,\dots,h\;.
$$
Then $n_k>1$ for all $k$.
Now we can state the result due to~\cite{S}, \cite{M} and~\cite{Eph}.
\begin{thm}\label{thm21}{\rm(Smith--Merle--Ephraim)}
Suppose that $f=0$ is a singular branch and $l=0$ a regular curve.
Let $\bar{b}_0,\dots,\bar{b}_h$ be the $(f,l)_0$-minimal system of
generators of the semigroup $S(f)$. Then with the notation introduced
above
$$
  \cQ(f,l)=\sum_{k=1}^h\left\{\Teisss{(n_k-1)\bar{b}_k}
                                       {(n_k-1)n_1\dots n_{k-1}}{18}{9}
                                       \right\}\;.
$$
\end{thm}
By convention the empty product which appears for $k=1$ is equal to 1.

The sequence of generators can be characterized in purely arithmetical
terms. Let us recall (see \cite{Br}, \cite{Zariski}, \cite{D}, \cite{GP2}).
\begin{thm}\label{thm22}
Let $\bar{b}_0,\bar{b}_1,\dots,\bar{b}_h$ be a sequence of strictly positive
integers. Then the following two conditions are equivalent.
\begin{itemize}
\item[\rm(I)] There is a singular branch $f=0$ and a regular curve $l=0$
such that $\bar{b}_0,\bar{b}_1,\dots,\bar{b}_h$ is the $(f,l)_0$-minimal
system of generators of the semigroup $S(f)$,
\item[\rm(II)] the sequence $\bar{b}_0,\bar{b}_1,\dots,\bar{b}_h$
satisfies the conditions:
\begin{itemize}
\item[\rm(Z$_1$)] the sequence $e_k=\GCD(\bar{b}_0,\dots,\bar{b}_k)$
$(k=0,1,\dots,h)$ is strictly decreasing and $e_h=1$.
\item[\rm(Z$_2$)] the sequence $e_{k-1}\bar{b}_k$ $(k=1,\dots,h)$ is
strictly increasing.
\end{itemize}
\end{itemize}
\end{thm}
\begin{ex}
{\rm For any integer $n\geq 0$ there is a branch $f=0$ with the semigroup
$\langle 6,8,27+6n\rangle$. By Theorem~\ref{thm21} we get
$$
  \cQ(f)=
2\,\left\{\Teisss{8}{1}{2}{1}\right\}+
3\,\left\{\Teisss{9+2n}{1}{6}{3}\right\}\;.
$$
} 
\end{ex}
Using Theorems~\ref{thm21} and~\ref{thm22} we get
\begin{cor}\label{cor24}\mbox{} Let $f=0$ be a singular branch. Then
\begin{itemize}
\item[\rm(1)] $\cQ(f,l)$ is a complete invariant of the pair $f=0,l=0$;
\item[\rm(2)] $\cQ(f)$ is a complete invariant of the branch $f=0$.
\end{itemize}
\end{cor}
\begin{thm}\label{thm25}{\rm\cite{GG}}\\
Let $f=0$ and $g=0$ be two reduced curves such that $\cQ(f)=\cQ(g)$.
Suppose that $f=0$ is an irreducible curve. Then $g=0$ is also irreducible.
\end{thm}
For every sequence $\bar{b}_0,\dots,\bar{b}_h$ satisfying conditions
(Z$_1$) and (Z$_2$) (in the sequel we call such a sequence~(Z)-sequence) we put
$$
  \cN(\bar{b}_0,\dots,\bar{b}_h)=
  \sum_{k=1}^h\left\{\Teisss{(n_k-1)\bar{b}_k}
  {(n_k-1)n_1\dots n_{k-1}}{18}{9}\right\}\;.
$$
and call $\cN(\bar{b}_0,\dots,\bar{b}_h)$ the Newton diagram associated with
the sequence $\bar{b}_0,\dots,\bar{b}_h$. Theorems~\ref{thm21}, \ref{thm25}
and Proposition~\ref{prop115} give rise to the following
\begin{cor}\label{cor26}{\rm(Irreducibility Criterion)}\\
Let $f=Y^n+a_1(X)Y^{n-1}+\dots+a_n(X)\in\bC\{X\}[Y]$ be a
distinguished polynomial of degree $n>1$ without multiple factors.
Then $f$ is irreducible if and only if the Newton diagram of the
discriminant $D(X,T)=\disc_Y(f(X,Y)-T)$ is equal to the Newton
diagram $\cN(\bar{b}_0,\dots,\bar{b}_h)$ associated with a
(Z)-sequence $\bar{b}_0=n,\bar{b}_1,\dots,\bar{b}_h$.
\end{cor}

\begin{ex}\label{ex37}
{\rm (see \cite{Kuo} and \cite{Ab}).
The following two examples are taken from \cite{Kuo}.

\vspace{2ex}\noindent \textbf{I.}
Let $f=(X^2-Y^3)^2-Y^7$. Then the curves $f=0$ and $Y=0$ are transverse
and $\cQ(f)=\cQ(f,Y)=\left\{\Teis{6}{1}\right\}+\left\{\Teis{14}{2}\right\}$.
To decide if $f$ is irreducible suppose that $\cQ(f)=\cN(\bar{b}_0,\dots,\bar{b}_h)$
for a (Z)-sequence $\bar{b}_0,\dots,\bar{b}_h$. Then $h=2$ since $\cQ(f)$ has two faces
and $\bar{b}_0=\ord1 f=4$. From condition
$$
\left\{\Teisss{(n_1-1)\bar{b}_1}{n_1-1}{12}{6}\right\} +
\left\{\Teisss{(n_2-1)\bar{b}_2}{(n_2-1)n_1}{12}{6}\right\} =
\left\{\Teisss{6}{1}{2}{1}\right\} +
\left\{\Teisss{14}{2}{3}{1.5}\right\}
$$
we get $\bar{b}_1=6$ and $\bar{b}_2=14$. A contradiction since
$\GCD(\bar{b}_0,\bar{b}_1,\bar{b}_2)=2$. Therefore $f$ is not irreducible.

\vspace{2ex}\noindent \textbf{II.} Let $f=(X^2-Y^3)^2-Y^5X$. The
curves $f=0$ and $Y=0$ are transverse and
$\cQ(f)=\cQ(f,Y)=\left\{\Teis{6}{1}\right\}+\left\{\Teis{13}{2}\right\}$
(see Example~\ref{e116}). It is easy to check that
$\cQ(f)=\cN(4,6,13)$ and that $4,6,13$ is a (Z)-sequence.
Therefore $f$ is irreducible with semigroup $S(f)=\langle
4,6,13\rangle$. }
\end{ex}

\section{Polar invariants in many branched case}

Let $\varphi,\psi\in\bC\{X,Y\}$ be irreducible power series.
The {\it contact coefficient\/} (in the sense of Hironaka) with
respect to a regular curve $l=0$ is the rational number
$$
  h(\varphi,\psi;l)=\frac{(\varphi,\psi)_0}{(l,\psi)_0}\;.
$$
If $l=0$ and $\psi=0$ are transverse then
$h(\varphi,\psi;l)=(\varphi,\psi)_0/\ord1\,\psi$
and we write $h(\varphi,\psi)$ instead of $h(\varphi,\psi;l)$.

Let $f=0$ be a reduced curve with $r>1$ branches. To describe the
contacts of $f_i=0$ with the branches $f_j=0$, $j\neq i$ let us
consider the following diagram
$$
  \cH_i(f,l)=\sum_{j=1}^r\left\{
               \Teisss{(f_i,f_j)_0}{(l,f_j)_0}{7}{3.5}
               \right\}
$$
and the set
$$
  H_i(f,l)=\left\{\frac{(f_i,f_j)_0}{(l,f_j)_0}:\;j\neq i\right\}\;.
$$
Note that the diagram $\cH_i(f,l)$ lies above horizontal axis and has
vertices $(0,(l,f)_0)$  and $((f_i,f/f_i)_0,(l,f_i)_0)$. The distance from
$\cH_i(f,l)$ to the horizontal axis is equal to $(l,f_i)_0$.
$$
\epsfbox{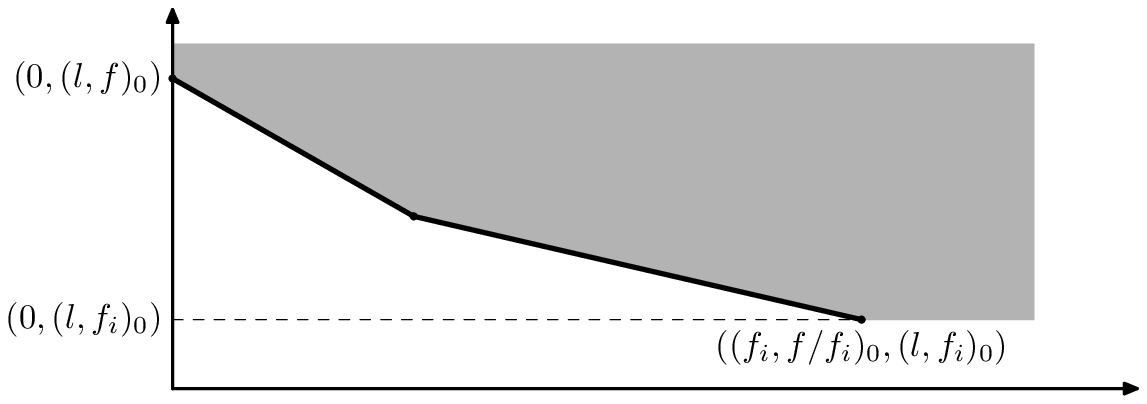}
$$
We omit the simple proof of the following
\begin{lm}\label{lm31}
The line with slope $-1/\tau$ ($\tau>0$) supporting $\cH_i(f,l)$
intersects the horizontal axis at the point
$$
  \left(\sum_{j=1}^r\inf\{(f_i,f_j)_0,\tau(l,f_j)_0\},0\right)\;.
$$
\end{lm}
Now let
$$
  q_i(\tau)=\frac{1}{(l,f_i)_0}
  \sum_{j=1}^r\inf\{(f_i,f_j)_0,\tau(l,f_j)_0\}
$$
for $\tau>0$ and $i=1,\dots,r$. According to Lemma~\ref{lm31} the function $q_i$ is determined
by the diagram $\cH_i(f,l)$ and has an obvious geometric interpretation.
The functions $q_i$ are piecewise linear, continuous and strictly increasing.
The following explicit formula for polar quotients of a many-branched curve
is due to \cite{GP3}.
\begin{thm}
Let $f=f_1\dots f_r$ be a reduced
power series with $r>1$ irreducible factors. Then
$$
  Q(f,l)=\bigcup q_i(Q(f_i,l)\cup H_i(f,l))\;.
$$
\end{thm}
We call the elements of
$q_i(Q(f_i,l)\cup H_i(f,l))$ polar invariants associated with the
branch $f_i=0$. A polar invariant can be associated with more than
one branch.

The polar invariants associated with the branch $f_i=0$ can be
interpreted in terms of the Newton diagram $\cH_i(f,l)$ and the
jacobian Newton polygon $\cQ(f_i,l)$ of the branch $f_i=0$. To
this end call a line supporting $\cH_i(f,l)$
\textit{distinguished} if it extends a face of $\cH_i(f,l)$ or is
parallel to a face of $\cQ(f_i,l)$. Then the polar invariants
associated to the branch $f_i=0$ are exactly the quotients of the
form $\frac{p}{d_i}$ where $(p,0)$ is the point of intersection of
a distinguished supporting line with the horizontal axis and
$d_i=(l,f_i)_0$ is the distance from $\cH_i(f,l)$ to this axis.

Let us calculate $\eta(f,l)=\sup Q(f,l)$. Using the fact that the
functions $q_i$ are increasing we get
\begin{thm}{\rm\cite{P2}}\\
$$
\eta(f,l)=
\begin{array}{c}r\\\max\\i=1\end{array}
\left\{
\max
\left\{
\eta(f_i,l),
\begin{array}{c}\mbox{}\\\max\\j\neq i\end{array}
\frac{(f_i,f_j)_0}{(l,f_j)_0}
\right\}+
\frac{1}{(l,f_i)_0}\sum_{j\neq i}(f_i,f_j)_0
\right\}\;.
$$
\end{thm}

For the applications of the above formula see \cite{GKP}.

If $f=0$ and $l=0$ are transverse then we write
$\cH_i(f)=\cH_i(f,l)$.
\begin{ex}
{\rm Even when $f=0$ and $g=0$ are curves with smooth branches the conditions
$\cH_i(f)=\cH_i(g)$ ($i=1,\dots,r$) do not imply the equisingularity of
$f=0$ and $g=0$. Let $f=f_1\dots f_{10}$ and $g=g_1\dots g_{10}$ where
$f_1=Y-X-X^2$,
$f_2=Y-X-2X^2$,
$f_3=Y-X-3X^2$,
$f_4=Y-2X-X^2$,
$f_5=Y-2X-2X^2$,
$f_6=Y-2X-3X^2$,
$f_7=Y-X$,
$f_8=Y-X-X^3$,
$f_9=Y-2X$,
$f_{10}=Y-2X-X^3$ and
$g_1=Y-X$,
$g_2=Y-X-X^2$,
$g_3=Y-X-2X^2$,
$g_4=Y-X-3X^2$,
$g_5=Y-X-4X^2$,
$g_6=Y-2X-2X^2$,
$g_7=Y-2X$,
$g_8=Y-2X-X^3$,
$g_9=Y-2X-X^2$,
$g_{10}=Y-2X-X^2-X^3$.
Then $\cH_i(f)=\cH_i(g)$ for $i=1,\dots,10$ but $f=0$ and $g=0$ are not
equisingular.
} 
\end{ex}

To construct a complete invariant of the pair $f=0$, $l=0$ the
notion of partial polar quotient introduced in \cite{E} is useful.
E.~Garc\'{\i}a Barroso characterized the type of equisingularity
of the curve by matrices of partial polar quotients (see
\cite{G}).

\section{Polar invariants and the Newton diagram}

We want to calculate the jacobian Newton polygon of a nondegenerate
singularity $f=0$ in terms of the Newton diagram $\Delta(f)$.
To formulate the result we need some notions. Let $f\in\bC\{X,Y\}$
be a nonzero power series without constant term.
The segment $S\in\cN_f$ is {\it principal\/} if $|S|_1=|S|_2$.
If a principal segment exists it is unique.
Put $\cN'_f=\cN_f\setminus\{\mbox{principal segment}\}$.
For every segment $S\in\cN'_f$ we put $m(S)=\min(|S|_1,|S|_2)-1$
if $1\leq|S|_1<|S|_2$ and $S$ has a vertex on the vertical axis
or if $1\leq|S|_2<|S|_1$ and $S$ has a vertex on the horizontal axis.
Moreover we let $m(S)=\min(|S|_1,|S|_2)$ for all remaining cases.

Let $\alpha/\alpha(S)+\beta/\beta(S)=1$ be the equation of the line
containing $S$. Obviously $\alpha(S),\beta(S)>0$ are rational numbers
and $\alpha(S)/\beta(S)=|S|_1/|S|_2$.

Recall that $t(f)$ is the number of tangents to $f=0$. If $f$ is
nondegenerate then $t(f)$ can be read from the Newton diagram $\Delta(f)$.
We have the following result due to \cite{LP} (see also \cite{LMP}).
\begin{thm}\label{thm41}
Suppose that $f$ is a nondegenerate singularity. Then
$$
  \cQ(f)=
  \left\{
  \Teisss{(\ord1\,f)(t(f)-1)}{(t(f)-1)}
  {16}{8}
  \right\}
  +
  \sum_{S\in\cN'_f}
  \left\{
  \Teisss{\max(\alpha(S),\beta(S))\,m(S)}{m(S)}
  {24}{12}
  \right\}\;.
$$
\end{thm}
We put by convention
$\{\Teis{0}{0}\}=\bR_+^2$ (the zero Newton diagram).
\begin{ex}\label{ex42}
{\rm Let $f=\sum c_{\alpha\beta}X^\alpha Y^\beta$ with summation
over all $(\alpha,\beta)\in\bN^2$ such that
$(\alpha/w_1)+(\beta/w_2)=1$ where $w_1,w_2\geq2$ are rational
numbers defines a reduced curve $f=0$. Then
$\eta_0(f)=\max(w_1,w_2)$ by Theorem~\ref{thm41}. On the other
hand $\mu_0(f)=(w_1-1)(w_2-1)$ by the Milnor-Orlik formula. Hence
the set of weights
$$
\{w_1,w_2\}=\left\{\frac{\mu_0(f)}{\eta_0(f)-1}+1,\,\eta_0(f)\right\}
$$
is an invariant of $f=0$.
} 
\end{ex}

\section{Application to pencils of plane curve singularities}

When studying the singularities at infinity of polynomials in two complex
variables of degree $N>1$ one considers the pencils of plane curves of the
form $f_t=f-tl^N$, $t\in\bC$ where $f,l\in\bC\{X,Y\}$ are coprime and a regular curve
$l=0$ is not a component of the local curve $f=0$ (see~\cite{Eph},
\cite{GP}, \cite{LMP}, \cite{P3}).
Let $U\subset\bC$ be a Zariski
open subset of $\bC$. We say that the pencil $(f_t:\,t\in U)$
is equisingular if the Milnor number $\mu_0(f_t)$ is constant
for $t\in U$. This means by $\mu$-constant theorem for
pencils~\cite{Casas2000} that for any $t_1,t_2\in U$ the curves
$f_{t_1}=0$ and $f_{t_2}=0$ are equisingular.
\begin{prop}\label{prop61} {\rm(\cite{Eph}, \cite{GP})}\\
Let $f=0$ be a reduced curve and $l=0$ a regular curve which is not
a branch of $f=0$. Let $N>0$ be an integer. Then
\begin{itemize}
\item[\rm(1)] the pencil $(f-tl^N:\,t\neq 0$) is equisingular if
and only if $N\not\in Q(f,l)$.
\item[\rm(2)] the pencil $(f-tl^N:\,t\in\bC$) is equisingular if
and only if $\eta(f,l)=\sup\,Q(f,l)<N$.
\end{itemize}
\end{prop}
Using the above proposition and a result of Ephraim~\cite{Eph} we get
the following
\begin{prop}\label{prop62}
Let $f=0$ be a singular branch, $l=0$ a regular one.
Let $(\bar{b}_0,\bar{b}_1,\dots,\bar{b}_h)_0$
be the $(f,l)_0$-minimal system of generators of the semigroup $S(f)$.
Then the following three conditions are equivalent
\begin{itemize}
\item[\rm(AM)] $e_{h-1}\bar{b}_h<(\bar{b}_0)^2$,
\item[\rm(I)] all series $f_t=f-tl^{\bar{b}_0}$, $t\in\bC$ are irreducible,
\item[\rm(E)] the pencil $(f_t=f-tl^{\bar{b}_0}:\,t\in\bC)$ is
equisingular.
\end{itemize}
\end{prop}

\noindent
\textbf{Proof.}
By Theorem~4.1 we have $\eta_0(f,l)=e_{h-1}\bar{b}_h/\bar{b}_0$.
Therefore (AM)~is equivalent to the inequality $\eta_0(f,l)<\bar{b}_0$ and
$\mbox{(AM)}\Leftrightarrow\mbox{(E)}$ follows from Proposition~\ref{prop61}(2).
Obviously $\mbox{(E)}\Rightarrow\mbox{(I)}$,
the implication $\mbox{(I)}\Rightarrow\mbox{(E)}$ is due to Ephraim
\cite{Eph}, Corollary~2.2.

\vspace{2ex}
Note that (AM) is the famous Abhyankar--Moh inequality
(see \cite{AM}, \cite{GP2}, \cite{Casas2000}).
For more applications of polar invariants to the singularities at infinity
we refer the reader to \cite{GP}, \cite{GP4}, \cite{P02}
and to the papers cited in these articles.

\vspace{2ex}\noindent
{\sc Department of Mathematics, Technical University,\\
AL. 1000 L PP 7, 25-314 Kielce, Poland\\
e-mail: }\verb+matjg@tu.kielce.pl ztpal@tu.kielce.pl matap@tu.kielce.pl+

\end{document}